\documentclass[12pt, a4paper]{article}
\usepackage{amsfonts,amssymb,amscd,amsmath,enumerate,verbatim,calc}

\newtheorem{thm2}{Theorem}[section]

\newtheorem{def2}{Definition}[section]
\begin{document}
\title{Continuity of Pseudo-Differential Operator $h_{\mu,a}$ Involving Hankel Translation and Hankel Convolution on Some Gevrey Spaces }
\author{Akhilesh Prasad and Manish Kumar \thanks{Department of Applied Mathematics, Indian School of Mines, Dhanbad-826004, India E-mail: apr$\_$bhu@yahoo.com, manish.math.bhu@gmail.com }}
\date{}
\maketitle
\begin{abstract}
The Pseudo-Differential Operator (p.d.o.) $h_{\mu,a}$ associated with the Bessel Operator involving the symbol $a(x,y)$
whose derivatives satisfy certain growth conditions depending on some increasing sequences is studied on certain Gevrey spaces. The p.d.o. $h_{\mu,a}$ on Hankel translation $\tau$ and Hankel convolution of Gevrey functions is continuous linear map into another Gevrey spaces.
\end{abstract}
{\bf Key Words:} Hankel transformation, Hankel translation, Hankel convolution, Pseudo-differential operator, Gevrey space.\\\\
{\bf MSC:} 46F05, 46F12.

\section{Introduction}
The pseudo-differential operator (p.d.o.) $h_{\mu,a}$ have applications in the study of boundary value problems on the half line. The p.d.o. $h_{\mu,a}$ was introduced by {\cite{pp}}, and its properties were investigated using Zemanian's theory of the Hankel transformation to certain space of Ultradistributations. Zemanian's theory was further extended by Lee \cite{le}, Pathak and Prasad \cite{pp}, for this purpose, the spaces $H_{\mu,a_k,A},H_\mu^{b_q,B}$ and $H_{\mu,a_k,A}^{b_q,B}$ of ultradifferentiable function were defined as follows. Similar spaces have been studied in {\cite{gs}} and {\cite{ro}}. Zemanian {\cite{ze}} introduced the function space $H_{\mu}$ consisting of all complex valued infinitely differentiable functions $\phi$ defined on $I=(0,\infty)$ satisfying
\begin{eqnarray}
\gamma_{m,k}^{\mu}(\phi)=\sup_{x\in I} \left|x^m (x^{-1}d/dx)^k x^{-\mu-1/2}\phi(x)\right| < \infty, \qquad \forall \ m,k \in {\bf N}_0.  \label{eq:1}
\end{eqnarray}
The pseudo-differential operator involving the symbol $a(x,y)$ is defined by
\begin{eqnarray}
(h_{\mu,a}\phi)(x)=\int_{0}^\infty (xy)^{1/2} J_{\mu}(xy) a(x,y) \widehat{\phi}(y) dy,  \qquad \mu {\geq-{1/2}} \label{eq:2}
\end{eqnarray}
where $\widehat{\phi}$ is the Hankel transformation defined by
\begin{eqnarray}
\widehat{\phi}(y)=(h_{\mu}\phi)(y)=\int_{0}^\infty (xy)^{1/2} J_{\mu}(xy) \phi(x) dx, \label{eq:3}
\end{eqnarray}
and $ J_{\mu} $ is the Bessel function of the first kind of order $ \mu $. We shall study the properties of the symbol   $a(x,y)$ in section 2. The following definitions and results will be needed in the sequel.

The space $L_{\mu}^p (\mu\geq{-{1/2}})$ is the set of all measurable functions $\phi$ on $I=(0,\infty)$ such that
\begin{eqnarray} 
\left\|\phi\right\|_\mu^p=\int_0^\infty \left|\phi(x)\right|^p x^{\mu+1/2} dx <\infty. \label{eq:4}
\end{eqnarray}
The Hankel translation of $\phi\in L_\mu^1(I)$ is defined by
\begin{eqnarray}
(\tau_z\phi)(w)=\int_{0}^\infty \phi(y) D_{\mu}(y,w,z) dy,\qquad \forall \ w,z\in I \label{eq:5}
\end{eqnarray}
where
\begin{eqnarray}
D_{\mu}(y,w,z)=\int_{0}^\infty t^{-\mu-1/2} j_\mu (yt) j_\mu (wt) j_\mu (zt) dt, \label{eq:6}
\end{eqnarray}
and
\begin{eqnarray}
j_\mu (wt)= (wt)^{1/2} J_\mu (wt). \label{eq:7}
\end{eqnarray}
The Hankel Convolution transform of two functions $\phi,\psi \in L_{\mu}^1(I)$ is defined by
\begin{eqnarray}
(\phi\#\psi)(z)=\int_0^\infty \phi (w) (\tau_z \psi) (w) dw, \qquad a.e. z\in I  \label{eq:8}
\end{eqnarray}
we shall also make use of the following results [{\cite {bm}},p. 285]
\begin{eqnarray}
h_{\mu}(\tau_z\phi)(u)=u^{-\mu-1/2} j_\mu(uz) (h_{\mu}\phi) (u), \qquad \forall \ u,z\in I \label{eq:9}
\end{eqnarray}
and
\begin{eqnarray}
h_{\mu}(\phi\#\psi)(u)=u^{-\mu-1/2}  (h_{\mu}\phi)(u) (h_{\mu}\psi)(u). \qquad \forall \ u\in I \label{eq:10}
\end{eqnarray}
We shall use the notation and terminology of {\cite{p,pp,ze}}. The differential operators $ N_\mu, M_\mu $ and $ S_\mu $ are defined by
\begin{eqnarray}
N_\mu=N_{\mu,x}=x^{\mu+1/2} (d/dx) x^{-\mu-1/2},  \label{eq:11} 
\end{eqnarray}
\begin{eqnarray}
M_\mu=M_{\mu,x}=x^{-\mu-1/2} (d/dx) x^{\mu+1/2}, \label{eq:12}
\end{eqnarray}
\begin{eqnarray}
 S_\mu=S_{\mu,x}=M_\mu N_\mu=d^2/dx^2+\frac{(1-4{\mu}^2)}{4x^2}. \label{eq:13}
\end{eqnarray}
We have the following relations for any $\phi \in H_{\mu} :$
\begin{eqnarray}
h_{\mu+1}(-x\phi)=N_\mu h_\mu \phi, \label{eq:14}
\end{eqnarray}
\begin{eqnarray}
h_{\mu+1}(N_\mu \phi)=-y h_{\mu}\phi, \label{eq:15}
\end{eqnarray}
\begin{eqnarray}
 h_{\mu}(S_{\mu}\phi)=-y^2 h_{\mu}\phi \label{eq:16}
\end{eqnarray}
and
\begin{eqnarray}
S_{\mu}^r \phi(x)=\sum_{j=0}^r b_j x^{2j+\mu+1/2} (d/dx)^{r+j} x^{-\mu-1/2}\phi(x), \label{eq:17}
\end{eqnarray}
where the $b_j$ are constants depending on $\mu$.

The following formula are given in [{\cite{ze}},pp.129,134] and [{\cite{p}},pp.240,242]
\begin{eqnarray}
(x^{-1}d/dx)^k (x^{-\mu-1/2}\psi\phi)=\sum_{\nu=0}^k\binom{k}{\nu}(x^{-1}d/dx)^\nu \psi (x^{-1}d/dx)^{k-\nu} (x^{-\mu-1/2}\phi) \label{eq:18}
\end{eqnarray}
\begin{eqnarray}
(x^{-1}d/dx)^k (x^{-\mu}J_\mu (x))=(-1)^k (x)^{-(\mu+k)}J_{\mu+k}(x) \label{eq:19}
\end{eqnarray}
\begin{eqnarray}
(x^{-1}d/dx)^k (x^{\mu}J_\mu (x))=(x)^{\mu-k}J_{\mu-k}(x). \label{eq:20} 
\end{eqnarray}
Let $\left\{a_k\right\}_{k \in{\bf N}_0} $ and $\left\{b_q\right\}_{q \in{\bf N}_0}$ be arbitrary sequences of positive numbers which satisfy the following conditions
\begin{eqnarray}
a_k^2\leq a_{k-1} a_{k+1},\qquad \qquad \forall \ k\geq 1 \label{eq:21}
\end{eqnarray}
\begin{eqnarray} 
b_q^2\leq b_{q-1} b_{q+1}\qquad \qquad \forall \ q\geq 1 \label{eq:22}
\end{eqnarray}
Immediate consequences of these inequalities are
\begin{eqnarray}
a_p a_k\leq a_{0} a_{p+k},\qquad \qquad \forall \ p,k=0,1,2,\dots \label{eq:23}
\end{eqnarray}
\begin{eqnarray}
b_p b_q\leq b_{0} b_{p+q}.\qquad \qquad \forall \ p,q=0,1,2,\dots \label{eq:24}
\end{eqnarray}
from inequality (~\ref{eq:21}) it can be proved that 
\begin{eqnarray}
(a_k/a_{k+1})\leq(a_{k-1}/a_k)\leq(a_{k-2}/a_{k-1})\dotsm \leq(a_{0}/a_{1}), \label{eq:25}
\end{eqnarray}
and
\begin{eqnarray}
a_{k-r}&=&(a_{k-r}/a_{k-r+1})\times(a_{k-r+1}/a_{k-r+2})\times \dotsm (a_{k+1}/a_k)\times a_k\nonumber\\
& \leq & {(a_{0}/a_{1})\times(a_{0}/a_{1})\times \dotsm (a_{0}/a_{1})} \times a_k;\nonumber
\end{eqnarray}
so that
\begin{eqnarray}
a_{k-r} & \leq & (a_{0}/a_{1})^r\times a_{k}. \label{eq:26}                                              \end{eqnarray}
Furthermore, assume that there are constants $R_1,R_2 > 0$ and $H_1,H_2 > 1$\\ 
such that
\begin{eqnarray}
a_p \leq R_1 H_1^p \min_{0\leq q\leq p} a_q a_{p-q}, \qquad \forall \ p,q \in {\bf N}_0  \label{eq:27}
\end{eqnarray}
\begin{eqnarray}
b_p \leq R_2 H_2^p \min_{0\leq q\leq p} b_q b_{p-q}, \qquad \forall \ p,q \in {\bf N}_0. \label{eq:28}
\end{eqnarray}
Let the constants $c_1,h_1,c_2,h_2,c$ and $h$ be such that for all $k,q\in{\bf N}_0,$ 
\begin{eqnarray}
a_{k+1}\leq c_1 h_1^k a_k, \label{eq:29}
\end{eqnarray}
\begin{eqnarray}
b_{q+1}\leq c_2 h_2^q b_q, \label{eq:30}
\end{eqnarray}
\begin{eqnarray}
b_{q+1}\geq c h^q b_q. \label{eq:31}
\end{eqnarray}
The conditions (~\ref{eq:29}) and (~\ref{eq:30}) may be replaced by the following stronger conditions whenever necessary
\begin{eqnarray}
a_{r+k}\leq L_1 R_1^{k+r} a_r a_k, \qquad \forall \ r,k\geq 0 \label{eq:32}
\end{eqnarray}
\begin{eqnarray}
b_{r+q}\leq L_2 R_2^{r+q} b_r b_q, \qquad \forall \ r,q\geq 0 \label{eq:33}
\end{eqnarray}
where $L_1,R_1,L_2$ and $R_2$ are positive constants.\\
The spaces of type $ H_{\mu}$, that is $ H_{\mu,a_k,A},H_{\mu}^{b_q,B}$ and $ H_{\mu,a_k,A}^{b_q,B}$ are defined as follows:
\begin{def2}
Let $\phi$ be infinitely differentiable function on I. Then $\phi\in H_{\mu,a_k,A}$ if and only if\\

$\left\|\phi\right\|_q^\mu=\sup_{k\in {\bf N}_0}\sup_{x\in I} \frac{\left|x^k(x^{-1}d/dx)^q x^{-\mu-1/2} \phi(x)\right|}{(A+\sigma)^k a_k} < \infty$\\\\
for every $q\in {\bf N}_0$ where $A$ is a certain positive constant depending on $\phi$ and $\sigma >0$ is arbitrary.
\end{def2}
\begin{def2}
The space $H_\mu^{b_q,B}$ is defined as follows: $\phi\in H_\mu^{b_q,B}$ if and only if\\

$\left\|\phi\right\|_k^\mu=\sup_{q\in {\bf N}_0}\sup_{x\in I} \frac{\left|x^k(x^{-1}d/dx)^q x^{-\mu-1/2} 
\phi(x)\right|}{(B+\rho)^q b_q} < \infty$\\\\
for every $k\in {\bf N}_0$ where $B$ is a positive constant depending on $\phi$ and $\rho >0$ is arbitrary.
\end{def2}
\begin{def2}
The function $\phi\in H_{\mu,a_k,A}^{b_q,B}$ if and only if\\

$\left\|\phi\right\|^\mu=\sup_{k,q\in {\bf N}_0}\sup_{x\in I} \frac {\left|x^k(x^{-1}d/dx)^q x^{-\mu-1/2} 
\phi(x)\right|}{(A+\sigma)^k a_k(B+\rho)^q b_q} < \infty$\\\\
where $\sigma $ and $\rho$ are as above and $A$ and $B$ are certain positive constants depending on $\phi$.\\
\end{def2}
The elements of the spaces $ H_{\mu,a_k,A}, H_{\mu}^{b_q,B}$ and $ H_{\mu,a_k,A}^{b_q,B}$ are called Ultradifferentiable functions and those of the corresponding dual spaces $(H_{\mu,a_k,A})', (H_{\mu}^{b_q,B})'$ and $(H_{\mu,a_k,A}^{b_q,B})'$ are called Ultradistributions.

From \cite{rk}, we have the following results:
\begin{thm2}
If $\left\{a_k\right\}$  satisfies (~\ref{eq:27}) and $\left\{b_q\right\}$ satisfies (~\ref{eq:28}) $ \forall \ k,q \in {{\bf N}_0}$, then for each fixed z, $0 < z < z_0 , \mu\geq{-1/2} $, the mapping $\phi \longmapsto \tau_{z} \phi $ is continuous from the spaces \\

(i) $H_{\mu,a_k,A}^{b_q,B}$ into $H_{\mu,a_k^3 b_k,A_3}^{a_q^2 b_q^2,B_4}$, where $A_3=A_1B_3 ({R}^{\otimes})^2,  B_3=R_1^2[B_1+(z_0 a_0/a_1)^2], A_1=AB (R^\ast)^2, B_4=A_1^2 (R^{\otimes})^6, R^{\otimes}=\max (1,R_1 R_2)$ and $ R_1, R_2$ are defined by (~\ref{eq:32}) and (~\ref{eq:33}), \\

and\\

(ii) $H_{\mu,a_k,A}$ into $H_{\mu,a_k^2,A_2}$, where $A_2=R_1^2[B_1+(z_0 a_0/a_1)^2], B_1=A^2(R^\ast)^6 $ and $ R^\ast=\max(1,R_1)$. 

We note that the Hankel translation cannot be defined on the whole of the space $H_{\mu}^{b_q,B}$; but it could be defined on a certain subspace $\widetilde{H}_{\mu}^{b_q,B}$ of $H_{\mu}^{b_q,B}$ in which the following condition is satisfied
\begin{eqnarray}
\sup_k {Q_{k+2q}^{\mu}}=Q^{\ast \mu}_q, \label{eq:34}
\end{eqnarray}
where $Q^{\ast \mu}_q$ are constants restraining the ${\phi}$'s in $H_{\mu}^{b_q,B}$. Then\\

(iii) $\widetilde{H}_{\mu}^{b_q,B}$ into $H_{\mu}^{b_q^2,B_2}$, where $B_2=B^2(R^\ast)^6, R^\ast=\max(1,R_1)$.
\end{thm2}

\section{Pseudo-Differential Operator Involving Hankel Translation on the Spaces of Type $H_\mu$}
This section investigates the p.d.o. involving Hankel translation $\tau$ on the spaces $H_{\mu,a_k,A}, H_\mu^{b_q,B}$ and $H_{\mu,a_k,A}^{b_q,B}$.

\begin{def2}
The symbol $a(x,y)$ is defined to be a complex valued function belonging to the space ${\bf C^\infty}(I\times I)$, such that its derivatives satisfy the growth condition
\begin{eqnarray}
\left| (x^{-1}d/dx)^\alpha (y^{-1}d/dy)^\nu a(x,y) \right| \leq L_m (C+\delta)^\alpha c_\alpha (D+\eta)^\nu d_\nu (1+y)^{m-\nu}  \label{eq:35}
\end{eqnarray}
for all $\alpha, \nu \in {\bf N}_0, \delta >0, \eta >0$ and $L_m >0$, where $m$ is a fixed real number, and $\left\{c_\alpha \right\}$ and $\left\{d_\nu \right\}$ are certain sequences of positive real numbers satisfying some of the conditions of type (~\ref{eq:21})-(~\ref{eq:31}). The set of all such symbols will be denoted  $S_{c_\alpha,d_\nu}^m$.\\

We have the following interesting results {\cite{rk}}.
\end{def2}
\begin{thm2}
If $\left\{a_k\right\}$ and $\left\{b_q\right\} \forall \ k,q \in {{\bf N}_0}$, satisfies (~\ref{eq:21}) and (~\ref{eq:22}) respectively then for each fixed $z, 0 < z < z_0$ and $\mu\geq{-1/2}$, the mapping $\phi \longmapsto h_{\mu}\tau_{z} \phi $ is linear and continuous from

(i) $H_{\mu,a_k,A}$ into $H_\mu^{a_q^2,B_3}$,  where $B_3=R_1^2\left[B_1+(z_0a_0/a_1)^2 \right]$ \\

(ii) $\widetilde{H}_\mu^{b_q,B}$ into $H_{\mu,b_k,B}$, \\

(iii) $H_{\mu,a_k,A}^{b_q,B}$ into $H_{\mu,a_k b_k,A_1}^{a_q^2,B_3}$, where $A_1=ABH_1^2, B_3=R_1^2\left[B_1+(z_0a_0/a_1)^2 \right]$ and $B_1$ as above.
\end{thm2} 
\begin{thm2}
Let $\left\{a_k\right\}, \left\{b_k\right\}, \left\{c_k\right\}$ and $\left\{d_k\right\}, k\in{\bf N}_0$, satisfy condition (~\ref{eq:26}), ${\mu}\geq {-1/2}$, and let the symbol $a(x,y)$ satisfy (~\ref{eq:35}) then the p.d.o. $\phi\longmapsto h_{\mu,a}\tau_z\phi$ is continuous linear mapping from $H_{\mu,a_k,A}^{b_q,B}$ into $H_{\mu,a^{\ast 3}_k b^{\ast}_k,A_6}^{a^{\ast 2}_q  b^{\ast 2}_q,B_6}$, where $a^\ast_k=\max_{k \in {\bf N}_0}(a_k,d_k), b^\ast_q=\max_{q \in {\bf N}_0}(b_q,c_q), A_6=((a^\ast_0/a^\ast_1)D+B_3)A_1, B_6= (b^\ast_0a^{\ast 2}_0/b^{\ast}_1a^{\ast 2}_1) C+H^6A_1^2 $ and $ B_3, A_1$ as above.
\end{thm2} 
{\bf Proof:}
Suppose that $\phi\in H_{\mu,a_k,A}^{b_q,B}$, then by Theorem (2.1)(iii) $(h_\mu \tau_z\phi)\in H_{\mu,a_k b_k,A_1}^{a_q^2,B_3}$, where $A_1=AB(R^\ast)^2$ and $B_3=R_1^2\left[B_1+(z_0 a_0/a_1)^2 \right], R^\ast=\max(1,R_1)$.\\

Now assume that
\begin{eqnarray*}
\Phi(x)&=&(h_{\mu,a} \tau_z\phi)(x)\\
&=&\int_0^\infty (xy)^{1/2} J_\mu(xy) a(x,y) (h_\mu \tau_z\phi)(y)dy. \\
\end{eqnarray*}
Using Zemanian's technique [\cite{ze},p. 144] we have
\begin{eqnarray}
N_\mu\Phi(x)&=&x^{\mu+1/2} (d/dx) x^{-\mu-1/2} \Phi(x)\nonumber\\
&=&x^{\mu+1+1/2}(x^{-1}d/dx)x^{-\mu-1/2}\Phi(x). \label{eq:36}
\end{eqnarray}
\begin{eqnarray*}
N_{\mu+1}N_\mu\Phi(x)&=&x^{\mu+1+1/2} (d/dx) x^{-(\mu+1)-1/2} N_\mu \Phi(x)\\ 
&=&x^{\mu+2+1/2}(x^{-1}d/dx) x^{-\mu-3/2} [x^{\mu+3/2}( x^{-1}d/dx) x^{-\mu-1/2} \Phi(x)]\\
&=&x^{\mu+2+1/2} (x^{-1}d/dx)^2 x^{-\mu-1/2} \Phi(x).
\end{eqnarray*}
Similary, using (~\ref{eq:18}), we have
\begin{eqnarray}
N_{\mu+q-1}\ldots N_\mu \Phi(x)&=& x^{\mu+q+1/2} (x^{-1}d/dx)^q  x^{-\mu-1/2} \Phi(x) \label{eq:37}\\
&=& x^{\mu+q+1/2} (x^{-1}d/dx)^q  x^{-\mu-1/2}\nonumber\\
&\times & \int_0^\infty (xy)^{1/2}J_\mu(xy) a(x,y) (h_\mu \tau_z\phi)(y)dy\nonumber\\
&=& x^{\mu+q+1/2} (x^{-1}d/dx)^q \int_0^\infty y^{1/2} x^{-\mu} J_\mu(xy) a(x,y) (h_\mu \tau_z\phi)(y)dy\nonumber\\
&=& x^{\mu+q+1/2} \int_0^\infty y^{1/2} \sum_{r=0}^q \binom{q}{r} (x^{-1}d/dx)^{q-r} x^{-\mu} J_\mu(xy)\nonumber\\
&\times & (x^{-1}d/dx)^r a(x,y) (h_\mu \tau_z\phi)(y)dy\nonumber\\
&=& x^{\mu+q+1/2} \int_0^\infty y^{1/2} \sum_{r=0}^q \binom{q}{r} (-y)^{q-r} x^{-\mu-q+r} J_{\mu+q-r}(xy)\nonumber\\
&\times & (x^{-1}d/dx)^r a(x,y) (h_\mu \tau_z\phi)(y)dy.\nonumber
\end{eqnarray}
Therefore,
\begin{eqnarray}
N_{\mu+q-1}\ldots N_\mu \Phi(x)&=& \sum_{r=0}^q\binom{q}{r}\int_0^\infty x^{r+1/2} y^{1/2}(x^{-1}d/dx)^r a(x,y)\nonumber\\
&\times & (h_\mu \tau_z\phi)(y)(-y)^{q-r}J_{\mu+q-r}(xy)dy \label{eq:38}\\
&=& \sum_{r=0}^q\binom{q}{r} x^r \int_0^\infty (xy)^{1/2}J_{\mu+q-r}(xy)\nonumber\\
&\times & [(x^{-1}d/dx)^r a(x,y)(-y)^{q-r} (h_\mu \tau_z\phi)(y)]dy.\nonumber\\
&=& \sum_{r=0}^q\binom{q}{r} x^r h_{\mu+q-r}\nonumber\\
&\times & [(x^{-1}d/dx)^r a(x,y)(-y)^{q-r} (h_\mu \tau_z\phi)(y)](x). \label{eq:39}
\end{eqnarray}
Using formula ${-x}h_\mu\phi=h_{\mu+1}(N_\mu\phi)$ in (~\ref{eq:39}), we get
\begin{eqnarray*}
(-x)N_{\mu+q-1}\ldots N_\mu \Phi(x)&=&\sum_{r=0}^q\binom{q}{r} x^r(-x) h_{\mu+q-r}[(x^{-1}d/dx)^r a(x,y)(-y)^{q-r} (h_\mu \tau_z\phi)(y)](x)\\
&=&\sum_{r=0}^q\binom{q}{r} x^r h_{\mu+q-r+1}N_{\mu+q-r}\\
&\times & [(x^{-1}d/dx)^r a(x,y)(-y)^{q-r} (h_\mu \tau_z\phi)(y)](x)\\
&=&\sum_{r=0}^q\binom{q}{r} x^r(-x)\int_0^\infty(xy)^{1/2}J_{\mu+q-r+1}(xy)N_{\mu+q-r}\\
&\times & [(x^{-1}d/dx)^r a(x,y)(-y)^{q-r} (h_\mu \tau_z\phi)(y)]dy\\
&=&\sum_{r=0}^q\binom{q}{r}\int_0^\infty x^r(xy)^{1/2}J_{\mu+q-r+1}(xy)y^{\mu+q-r+1/2}(d/dy)\\
&\times & y^{-\mu-q+r-1/2}[(x^{-1}d/dx)^r a(x,y)(-y)^{q-r} (h_\mu \tau_z\phi)(y)]dy
\end{eqnarray*}
\begin{eqnarray}
&=& \sum_{r=0}^q\binom{q}{r}(-1)^{q-r}\int_0^\infty x^{r+1/2}y^{\mu+q-r+2} (y^{-1}d/dy)\nonumber\\
&\times & [y^{-\mu-1/2}(h_\mu \tau_z\phi)(y)(x^{-1}d/dx)^r a(x,y)J_{\mu+q-r+1}(xy)]dy \label{eq:40}
\end{eqnarray}
\begin{eqnarray}
&=& \sum_{r=0}^q\binom{q}{r}(-1)^{q-r}x^r\int_0^\infty(xy)^{1/2}J_{\mu+q-r+1}(xy)\nonumber\\
&\times & [y^{\mu+q-r+1+1/2}(y^{-1}d/dy)\big\{ y^{-\mu-1/2}(h_\mu \tau_z\phi)(y)(x^{-1}d/dx)^r a(x,y)\big\}]dy \nonumber\\
&=&\sum_{r=0}^q\binom{q}{r}(-1)^{q-r}x^r h_{\mu+q-r+1}[y^{\mu+q-r+1+1/2}\nonumber\\
&\times & (y^{-1}d/dy)\big\{ y^{-\mu-1/2}(h_\mu \tau_z\phi)(y)(x^{-1}d/dx)^r a(x,y)\big\}]. \label{eq:41}
\end{eqnarray}
Now, from (~\ref{eq:41}) again using result ${-x}h_\mu\phi=h_{\mu+1}(N_\mu\phi)$, we get
\begin{eqnarray*}
(-x)^2 (N_{\mu+q-1}\ldots N_\mu \Phi(x))&=&\sum_{r=0}^q\binom{q}{r}(-1)^{q-r}x^r (-x) h_{\mu+q-r+1}\\
&\times &[y^{\mu+q-r+1+1/2}(y^{-1}d/dy)\big\{ y^{-\mu-1/2}(h_\mu \tau_z\phi)(y)(x^{-1}d/dx)^r a(x,y)\big\}]\\
&=&\sum_{r=0}^q\binom{q}{r}(-1)^{q-r}x^r h_{\mu+q-r+2}N_{\mu+q-r+1}\\
&\times &[y^{\mu+q-r+1+1/2}(y^{-1}d/dy)\big\{ y^{-\mu-1/2}(h_\mu \tau_z\phi)(y)(x^{-1}d/dx)^r a(x,y)\big\}]\\
&=&\sum_{r=0}^q\binom{q}{r}(-1)^{q-r}x^r\int_0^\infty(xy)^{1/2}J_{\mu+q-r+2}(xy)N_{\mu+q-r+1}\\
&\times & [y^{\mu+q-r+1+1/2}(y^{-1}d/dy)\big\{ y^{-\mu-1/2}(h_\mu \tau_z\phi)(y)(x^{-1}d/dx)^r a(x,y)\big\}]\\
&=&\sum_{r=0}^q\binom{q}{r}(-1)^{q-r}\int_0^\infty x^{r+1/2} y^{\mu+q-r+2+1}\\
&\times &[(y^{-1}d/dy)^2\big\{ y^{-\mu-1/2}(h_\mu \tau_z\phi)(y)(x^{-1}d/dx)^r a(x,y)\big\}]\\
&\times & J_{\mu+q-r+2}(xy)dy.
\end{eqnarray*}
In general, we have
\begin{eqnarray*}
(-x)^k (N_{\mu+q-1}\ldots N_\mu \Phi(x))&=&\sum_{r=0}^q\binom{q}{r}(-1)^{q-r}\int_0^\infty x^{r+1/2} y^{\mu+q-r+k+1}\\
&\times &[(y^{-1}d/dy)^k\big\{ y^{-\mu-1/2}(h_\mu \tau_z\phi)(y)(x^{-1}d/dx)^r a(x,y)\big\}]\\
&\times & J_{\mu+q-r+k}(xy)dy.
\end{eqnarray*}
\begin{eqnarray}
&=&\sum_{r=0}^q\binom{q}{r}(-1)^{q-r}\int_0^\infty x^{r+1/2}y^{\mu+q-r+k+1}\nonumber\\
&\times & \sum_{\nu=0}^k\binom{k}{\nu}(y^{-1}d/dy)^\nu(x^{-1}d/dx)^r \nonumber\\
&\times & a(x,y)(y^{-1}d/dy)^{k-\nu}y^{-\mu-1/2}(h_\mu \tau_z\phi)(y)J_{\mu+q-r+k}(xy)dy. \label{eq:42}
\end{eqnarray}
Now, from (~\ref{eq:37}), we know that
\begin{eqnarray}
N_{\mu+q-1}\ldots N_\mu \Phi(x) &=& x^{\mu+q+1/2} (x^{-1}d/dx)^q x^{-\mu-1/2} \Phi(x). \label{eq:43}
\end{eqnarray}
Multiplying both sides in (~\ref{eq:43}) by $(-x)^k$, we get
\begin{eqnarray}
(-x)^k (N_{\mu+q-1}\ldots N_\mu \Phi(x))=(-1)^k x^{\mu+k+q+1/2} (x^{-1}d/dx)^q  x^{-\mu-1/2} \Phi(x). \label{eq:44}
\end{eqnarray}
Comparing equations (~\ref{eq:42}) and (~\ref{eq:44}), we have
\begin{eqnarray*}
(-1)^k x^{\mu+k+q+1/2} (x^{-1}d/dx)^q  x^{-\mu-1/2} \Phi(x)&=&\sum_{r=0}^q\binom{q}{r}(-1)^{q-r}\int_0^\infty x^{r+1/2}y^{\mu+q-r+k+1}\\
&\times & \sum_{\nu=0}^k\binom{k}{\nu}(y^{-1}d/dy)^\nu(x^{-1}d/dx)^r a(x,y)\\
&\times & (y^{-1}d/dy)^{k-\nu} y^{-\mu-1/2}(h_\mu \tau_z\phi)(y)J_{\mu+q-r+k}(xy)dy.
\end{eqnarray*}
Therefore,
\begin{eqnarray*}
(-1)^k x^{k} (x^{-1}d/dx)^q  x^{-\mu-1/2}\Phi(x)&=& \sum_{r=0}^q\binom{q}{r}(-1)^{q-r}\int_0^\infty x^{-(\mu+q-r)}y^{\mu+q-r+k+1}\\
&\times &\sum_{\nu=0}^k\binom{k}{\nu}(y^{-1}d/dy)^\nu(x^{-1}d/dx)^r a(x,y)\\
& \times & (y^{-1}d/dy)^{k-\nu}y^{-\mu-1/2}(h_\mu \tau_z\phi)(y)J_{\mu+q-r+k}(xy)dy.
\end{eqnarray*}
Thus
\begin{eqnarray*}
\left|x^{k} (x^{-1}d/dx)^q  x^{-\mu-1/2}\Phi(x)\right| &\leq & \sum_{r=0}^q\binom{q}{r} \int_0^\infty y^{2({\mu+q-r})+{k+1}}\sum_{\nu=0}^k\binom{k}{\nu}\\
&\times & \left|(x^{-1}d/dx)^r(y^{-1}d/dy)^\nu a(x,y)\right| \\
& \times & \left|(y^{-1}d/dy)^{k-\nu} y^{-\mu-1/2}(h_\mu\tau_z\phi)(y)\right|\\
&\times & \left|(xy)^{-(\mu+q-r)}J_{\mu+q-r+k}(xy)\right|dy.
\end{eqnarray*}
Using inequality (~\ref{eq:35}), the right-hand side assumes the form
\begin{eqnarray*}
\sum_{r=0}^q\binom{q}{r}\int_0^\infty y^{2({\mu+q-r})+{k+1}} \sum_{\nu=0}^k\binom{k}{\nu}L_m (C+\delta)^r c_r (D+\eta)^\nu d_\nu (1+y)^{m-\nu}\\
\times \left|(y^{-1}d/dy)^{k-\nu}y^{-\mu-1/2}(h_\mu \tau_z\phi)(y)\right|\frac{2^{-(\mu+q-r)}E}{\Gamma{(\mu+q-r+1})}dy.
\end{eqnarray*}
If we assume that $p$ is a positive integer such that $p\geq m$ and $ s > 2{\mu}+1 $, then the last term can be estimated by
\begin{eqnarray*}
\left|x^{k} (x^{-1}d/dx)^q  x^{-\mu-1/2}\Phi(x)\right| &\leq & \frac{2^{-\mu}EL_m}{\Gamma{(\mu+1)}} \sum_{r=0}^q\sum_{\nu=0}^k\binom{q}{r}\binom{k}{\nu} (C+\delta)^r c_r (D+\eta)^\nu d_\nu \\
&\times & \int_0^\infty y^{2({\mu+q-r})+{k+1}}(1+y)^{m-\nu+s}\\
&\times & \left|(y^{-1}d/dy)^{k-\nu}y^{-\mu-1/2}(h_\mu \tau_z\phi)(y)\right|\frac{dy}{(1+y)^s}\\
&\leq& \frac{2^{-\mu}EL_m}{\Gamma{\mu+1}} \sum_{r=0}^q\sum_{\nu=0}^k\binom{q}{r}\binom{k}{\nu} (C+\delta)^r c_r (D+\eta)^\nu d_\nu \\
&\times & \sup_{y\in I}[(1+y)^{p+s}y^{2(q-r)+k} \left|(y^{-1}d/dy)^{k-\nu}y^{-\mu-1/2}(h_\mu \tau_z\phi)(y)\right|]\\
&\times & \int_0^\infty \frac{y^{2(\mu+1/2)}}{(1+y)^s}dy.\\
\end{eqnarray*}
\begin{eqnarray*}
\left|x^{k} (x^{-1}d/dx)^q  x^{-\mu-1/2}\Phi(x)\right| &\leq & \frac{2^{-\mu}EL_m}{\Gamma{(\mu+1})} \sum_{r=0}^q\sum_{\nu=0}^k\binom{q}{r}\binom{k}{\nu} (C+\delta)^r c_r (D+\eta)^\nu d_\nu \\
&\times & \sum_{n=0}^{p+s}\binom{p+s}{n} \sup_{y\in I} [y^n y^{2(q-r)+k}\\
&\times & \left|(y^{-1}d/dy)^{k-\nu}y^{-\mu-1/2}(h_\mu \tau_z\phi)(y)\right|]\\
\end{eqnarray*}
\begin{eqnarray}
&\leq &\frac{2^{-\mu}EL_m}{\Gamma{(\mu+1)}} \sum_{r=0}^q\sum_{\nu=0}^k\sum_{n=0}^{p+s}\binom{q}{r}\binom{k}{\nu}\binom{p+s}{n} \nonumber\\
&\times &  (C+\delta)^r c_r (D+\eta)^\nu d_\nu \nonumber\\
& \times & \sup_{y\in I} \left|y^{n+2(q-r)+k}(y^{-1}d/dy)^{k-\nu}y^{-\mu-1/2}(h_\mu \tau_z\phi)(y)\right|. \label{eq:45}
\end{eqnarray}
Using Theorem 2.1 (iii) in (~\ref{eq:45}), we have
\begin{eqnarray*}
\left|x^{k} (x^{-1}d/dx)^q  x^{-\mu-1/2}\Phi(x)\right| &\leq &\frac{2^{-\mu}EL_m}{\Gamma{(\mu+1)}} \sum_{r=0}^q\sum_{\nu=0}^k\sum_{n=0}^{p+s}\binom{q}{r}\binom{k}{\nu}\binom{p+s}{n} \nonumber\\
&\times &  (C+\delta)^r c_r (D+\eta)^\nu d_\nu \left\|(h_\mu \tau_z\phi)(y)\right\|^\mu (A_1+\sigma)^{n+2(q-r)+k} \\
&\times & a_{n+2(q-r)+k}b_{n+2(q-r)+k}(B_3+\rho)^{k-\nu} a_{k-\nu}^2.
\end{eqnarray*}
Using inequality (~\ref{eq:27}) and (~\ref{eq:28}), we have
\begin{eqnarray*}
\left|x^{k} (x^{-1}d/dx)^q  x^{-\mu-1/2}\Phi(x)\right| &\leq &\frac{2^{-\mu}EL_m}{\Gamma{(\mu+1)}} \sum_{r=0}^q\sum_{\nu=0}^k\sum_{n=0}^{p+s}\binom{q}{r}\binom{k}{\nu}\binom{p+s}{n} 
(C+\delta)^r c_r (D+\eta)^\nu d_\nu \\
&\times & \left\|(h_\mu \tau_z\phi)\right\|^\mu (A_1+\sigma)^{n+2(q-r)+k} (B_3+\rho)^{k-\nu} R_{1}^3 H_{1}^{n+6(q-r)+2k}\\
&\times & a_n a_{q-r}^2 a_k R_{2}^3 H_{2}^{n+6(q-r)+2k} b_n b_{q-r}^2 b_k a_{k-\nu}^2.
\end{eqnarray*}
Applying inequality (~\ref{eq:23}) and (~\ref{eq:24}), we have
\begin{eqnarray*}
\left|x^{k} (x^{-1}d/dx)^q  x^{-\mu-1/2}\Phi(x)\right| &\leq &\frac{2^{-\mu}EL_m}{\Gamma{(\mu+1)}}(R_1 R_2)^3 \sum_{r=0}^q\sum_{\nu=0}^k\sum_{n=0}^{p+s}\binom{q}{r}\binom{k}{\nu}\binom{p+s}{n}
(C+\delta)^r c_r \\
&\times & (D+\eta)^\nu d_\nu H^{n+6(q-r)+2k} (A_1+\sigma)^{n+2(q-r)+k} a^\ast_n a^\ast_k a^{\ast 2}_{q-r}\\
&\times & (a^\ast_{k-\nu} a^\ast_\nu)a^\ast_{k-\nu}b^\ast_n b^{\ast}_{q-r}(b^{\ast}_{q-r}b^\ast_r)(B_3+\rho)^{k-\nu}b^\ast_k\left\|(h_\mu \tau_z\phi)\right\|^\mu
\end{eqnarray*}
where $b^\ast_q=\max_{q\in {\bf N}_0}(b_q,c_q),a^\ast_k=\max_{k\in {\bf N}_0}(a_k,d_k),H=H_1 H_2$ and $R_3=R_1 R_2$.\\

Finally
\begin{eqnarray*}
\left|x^{k} (x^{-1}d/dx)^qx^{-\mu-1/2}\Phi(x)\right| &\leq& \frac{2^{-\mu}EL_m}{\Gamma{(\mu+1)}} (R_3)^3\sum_{r=0}^q\sum_{\nu=0}^k\sum_{n=0}^{p+s}\binom{q}{r}\binom{k}{\nu}\binom{p+s}{n}(C+\delta)^r\\ 
&\times & H^{n+6(q-r)+2k} (D+\eta)^\nu(A_1+\sigma)^{n+2(q-r)+k}(B_3+\rho)^{k-\nu}a^\ast_n b^\ast_n b^\ast_k\\
&\times & (a^\ast_0/a^\ast_1)^{2r} a^{\ast 2}_q(a^\ast_0/a^\ast_1)^{\nu}a^\ast_0 a^\ast_k(b^\ast_0/b^\ast_1)^{r} b^\ast_q b^\ast_0 b^\ast_q\left\|(h_\mu \tau_z\phi)\right\|^\mu
\end{eqnarray*}
\begin{eqnarray*}
&\leq & (R_3)^3\sum_{r=0}^q\binom{q}{r}[(b^\ast_0 a^{\ast 2}_0/b^\ast_1 a^{\ast 2}_1)(C+\delta)]^r[H^6(A_1+\sigma)^2]^{q-r}\\
&\times & a^{\ast 2}_q b^{\ast 2}_q \sum_{\nu=0}^k \binom{k}{\nu}(B_3+\rho)^{k-\nu}[(a^\ast_0/a^\ast_1)(D+\eta)]^\nu (A_1+\sigma)^k a^{\ast 3}_k b^\ast_k\left\|(h_\mu \tau_z\phi)\right\|^\mu \\
&\leq & R (A_6+\sigma_1)^k a^{\ast 3}_k b^\ast_k (B_6+\rho_1)^q a^{\ast 2}_q b^{\ast 2}_q \left\|(h_\mu \tau_z\phi)\right\|^\mu
\end{eqnarray*}
where $A_6=((a^\ast_0/a^\ast_1)D+B_3)A_1, B_6= (b^\ast_0a^{\ast 2}_0/b^{\ast}_1a^{\ast 2}_1) C+H^6A_1^2$ and $R$ is a constant. Therefore
\begin{eqnarray*}
\left\|\Phi\right\|^\mu=\sup_{k,q\in {\bf N}_0}\sup_{x\in I}\frac{\left|x^{k} (x^{-1}d/dx)^qx^{-\mu-1/2}\Phi(x)\right|}{(A_6+\sigma_1)^k a^{\ast 3}_k b^\ast_k (B_6+\rho_1)^q a^{\ast 2}_q b^{\ast 2}_q}\leq R \left\|(h_\mu \tau_z\phi)\right\|^\mu.
\end{eqnarray*}
Hence
$\Phi(x)\in H_{\mu,a^{\ast 3}_k b^{\ast}_k,A_6}^{a^{\ast 2}_q  b^{\ast 2}_q,B_6}$.

{\bf REMARK:} In Theorem 2.2 we may choose $\phi \in H_{\mu,a_k,A}(or \widetilde{H}_\mu^{b_q,B})$ then the p.d.o. $h_{\mu,a}$ on Hankel translation $(h_{\mu,a}(\tau_z \phi)) \in \widetilde{H}_{\mu,a^{\ast 2}_k,A'_6}(or \widetilde{H}_{\mu}^{b^{\ast 2}_q,B'_6})$ where $A'_6=(a^\ast_0/a^\ast_1)D+B_3$ and $ B'_6=(b^{\ast}_0/b^{\ast}_1)+H_2^6 B^2$.

\section{Pseudo-Differential Operator Involving Hankel Convolution on the Spaces of Type $ H_{\mu} $ }
In this section we investigate the p.d.o. involving Hankel convolution transform of $ \phi\#\psi $ on the spaces $ H_{\mu,a_k,A},H_{\mu}^{b_q,B}$ and $H_{\mu,a_k,A}^{b_q,B}$.
We have the following interesting results {\cite{rk}}.

\begin{thm2}
If $\left\{a_k\right\}$ and $\left\{b_q\right\} \forall \ k,q \in {{\bf N}_0}$ satisfies (~\ref{eq:21}) and (~\ref{eq:22}) respectively then for $\mu \geq {-1/2}$, the mapping $(\phi,\psi)\longmapsto (\phi\#\psi)$ is linear and continuous from the spaces\\

(i) $H_{\mu,a_k,A}\times H_{\mu,a_k,A}$ into $\widetilde{H}_{\mu,a_k^2,B_1}$, where $B_1=A^2(R^\ast)^6$\\

(ii)$\widetilde{H}_{\mu}^{b_q,B} \times \widetilde{H}_{\mu}^{b_q,B}$ into $H_{\mu}^{b_q^2,B_2}$, where $B_2=B^2(R^\ast)^6$\\

(iii) $H_{\mu,a_k,A}^{b_q,B}\times H_{\mu,a_k,A}^{b_q,B}$ into  $H_{\mu,a_k^3 b_k,A_4}^{a_q^2 b_q^2,B_5}$, where $A_4=(R^\otimes)^2 A_1 B_1, B_5=(R^\otimes)^6 A_1^2,$ and $ R^\otimes=\max(1,R_1 R_2)$.
\end{thm2}

\begin{thm2}
If $\left\{a_k\right\}$ and $\left\{b_q\right\} \forall \ k,q \in {{\bf N}_0}$ satisfies (~\ref{eq:21}) and (~\ref{eq:22}) respectively then for $\mu \geq {-1/2}$, the mapping $(\phi,\psi)\longmapsto h_{\mu}(\phi\#\psi)$ is continuous linear mapping from\\

(i) $H_{\mu,a_k,A}\times H_{\mu,a_k,A}$ into $H_{\mu}^{a_q^2,B_1}$,\\
where $B_1=A^2(R^\ast)^6$ and $ R^\ast=\max(1,R_1)$ \\

(ii)$\widetilde{H}_\mu^{b_q,B} \times \widetilde{H}_\mu^{b_q,B}$ into $H_{\mu,b_k,B}$, and \\

(iii) $H_{\mu,a_k,A}^{b_q,B}\times H_{\mu,a_k,A}^{b_q,B}$ into $H_{\mu,a_k b_k,A_1}^{a_q^2,B_1}$, where $A_1=AB(R^\ast)^2, B_1=A^2(R^{\ast})^6$ and $R^\ast=\max(1,R_1)$.
\end{thm2} 
\begin{thm2}
Let $ \left\{a_k\right\}, \left\{b_k\right\}, \left\{c_k\right\}$ and $ \left\{d_k\right\}, k\in {\bf N}_0 $, satisfy condition (~\ref{eq:26}), $ \mu \geq {-1/2} $ and the symbol $a(x,y)$ satisfy (~\ref{eq:35}) then the p.d.o $ (\phi,\psi)\longmapsto h_{\mu,a}(\phi\#\psi)$ is continuous linear mapping from $ H_{\mu,a_k,A}^{b_q,B}\times H_{\mu,a_k,A}^{b_q,B}$ into $ H_{\mu,a^{\ast 3}_k b^{\ast}_k,A_7}^{a^{\ast 2}_q  b^{\ast 2}_q,B_6}$, where $A_7=((a^\ast_0/a^\ast_1) D+B_1)A_1$ and $ B_6= (b^\ast_0 a^{\ast 2}_0/b^{\ast}_1a^{\ast 2}_1) C+H^6 A_1^2$.
\end{thm2}
{\bf Proof:}
Let $(\phi,\psi)\in H_{\mu,a_k,A}^{b_q,B}\times H_{\mu,a_k,A}^{b_q,B}$ then by, Theorem 3.2 (iii),   $h_\mu (\phi\#\psi)\in H_{\mu,a_k b_k,A_1}^{a_q^2,B_1}$.

 Now assume that $\Psi(x)=h_{\mu,a} (\phi\#\psi)(x)$, 
from inequality (~\ref{eq:45}), we have
\begin{eqnarray*}
\left|x^{k} (x^{-1}d/dx)^q  x^{-\mu-1/2}\Psi(x)\right|&\leq &\frac{2^{-\mu}EL_m}{\Gamma{(\mu+1)}} \sum_{r=0}^q\sum_{\nu=0}^k\sum_{n=0}^{p+s}\binom{q}{r}\binom{k}{\nu}\binom{p+s}{n} \\
&\times &  (C+\delta)^r c_r (D+\eta)^\nu d_\nu \\
& \times & \sup_{y\in I} \left|y^{n+2(q-r)+k}(y^{-1}d/dy)^{k-\nu}y^{-\mu-1/2}h_{\mu} (\phi\#\psi)(y)\right|.
\end{eqnarray*} 
Using Theorem 3.2 (iii), the right-hand side assumes the form
\begin{eqnarray*}
\left|x^{k} (x^{-1}d/dx)^q  x^{-\mu-1/2}\Psi(x)\right| &\leq &\frac{2^{-\mu}EL_m}{\Gamma{(\mu+1)}} \sum_{r=0}^q\sum_{\nu=0}^k\sum_{n=0}^{p+s}\binom{q}{r}\binom{k}{\nu}\binom{p+s}{n}\\
&\times &  (C+\delta)^r c_r (D+\eta)^\nu d_\nu \left\|h_{\mu} (\phi\#\psi)\right\|^\mu (A_1+\sigma)^{n+2(q-r)+k} \\
&\times & a_{n+2(q-r)+k} b_{n+2(q-r)+k} (B_1+\rho)^{k-\nu} a_{k-\nu}^2.
\end{eqnarray*}
Using the inequalities (~\ref{eq:27}) and (~\ref{eq:28}) we can bound this expression by
\begin{eqnarray*}
\frac {2^{-\mu}EL_m}{\Gamma{(\mu+1)}} \sum_{r=0}^q\sum_{\nu=0}^k\sum_{n=0}^{p+s}\binom{q}{r}\binom{k}{\nu}\binom{p+s}{n} 
(C+\delta)^r c_r (D+\eta)^\nu d_\nu \left\|h_{\mu} (\phi\#\psi)\right\|^\mu(A_1+\sigma)^{n+2(q-r)+k} \\
\times(B_1+\rho)^{k-\nu} R_{1}^3 H_{1}^{n+6(q-r)+2k} a_n a_{q-r}^2 a_k R_{2}^3 H_{2}^{n+6(q-r)+2k} b_n b_{q-r}^2 b_k a_{k-\nu}^2.
\end{eqnarray*} 
Applying inequality (~\ref{eq:23}) and (~\ref{eq:24}), we have
\begin{eqnarray*}
\left|x^{k} (x^{-1}d/dx)^q  x^{-\mu-1/2}\Psi(x)\right| &\leq &\frac{2^{-\mu}EL_m}{\Gamma{(\mu+1)}}(R_1 R_2)^3 \sum_{r=0}^q\sum_{\nu=0}^k\sum_{n=0}^{p+s}\binom{q}{r}\binom{k}{\nu}\binom{p+s}{n} 
(C+\delta)^r c_r \\
&\times & (D+\eta)^\nu d_\nu H^{n+6(q-r)+2k} (A_1+\sigma)^{n+2(q-r)+k} a^\ast_n a^\ast_k a^{\ast 2}_{q-r}\\
&\times & (a^\ast_{k-\nu} a^\ast_\nu)a^\ast_{k-\nu}b^\ast_n b^{\ast}_{q-r}(b^{\ast}_{q-r}b^\ast_r)(B_1+\rho)^{k-\nu}b^\ast_k\left\|h_{\mu} (\phi\#\psi)\right\|^\mu,
\end{eqnarray*}
where $b^\ast_q=\max_{q\in {\bf N}_0}(b_q,c_q),a^\ast_k=\max_{k\in {\bf N}_0}(a_k,d_k),H=H_1 H_2$ and $R_3=R_1 R_2$.\\

Finally
\begin{eqnarray*}
\left|x^{k} (x^{-1}d/dx)^qx^{-\mu-1/2}\Psi(x)\right| &\leq& \frac{2^{-\mu}EL_m}{\Gamma{(\mu+1)}} (R_3)^3\sum_{r=0}^q\sum_{\nu=0}^k\sum_{n=0}^{p+s}\binom{q}{r}\binom{k}{\nu}\binom{p+s}{n}(C+\delta)^r\\ 
&\times & H^{n+6(q-r)+2k} (D+\eta)^\nu(A_1+\sigma)^{n+2(q-r)+k}(B_1+\rho)^{k-\nu}a^\ast_n b^\ast_n b^\ast_k\\
&\times & (a^\ast_0/a^\ast_1)^{2r} a^{\ast 2}_q(a^\ast_0/a^\ast_1)^{\nu}a^\ast_0 a^\ast_k(b^\ast_0/b^\ast_1)^{r} b^\ast_q b^\ast_0 b^\ast_q\left\|h_{\mu} (\phi\#\psi)\right\|^\mu
\end{eqnarray*}
\begin{eqnarray*}
&\leq & (R_3)^3\sum_{r=0}^q\binom{q}{r}[(b^\ast_0 a^{\ast 2}_0/b^\ast_1 a^{\ast 2}_1)(C+\delta)]^r[H^6(A_1+\sigma)^2]^{q-r}\\
&\times & a^{\ast 2}_q b^{\ast 2}_q \sum_{\nu=0}^k \binom{k}{\nu}(B_1+\rho)^{k-\nu}[(a^\ast_0/a^\ast_1)(D+\eta)]^\nu (A_1+\sigma)^k a^{\ast 3}_k b^\ast_k\left\|h_{\mu} (\phi\#\psi)\right\|^\mu \\
&\leq & R (A_7+\sigma_1)^k a^{\ast 3}_k b^\ast_k (B_6+\rho_1)^q a^{\ast 2}_q b^{\ast 2}_q \left\|h_{\mu} (\phi\#\psi)\right\|^\mu,
\end{eqnarray*}
where $ A_7=((a^\ast_0/a^\ast_1)D+B_1)A_1, B_6=( b^\ast_0 a^{\ast 2}_0/ b^\ast_1 a^{\ast 2}_1) C+H^6 A_1^2$ and $R$ is a constant. Therefore
\begin{eqnarray*}
\left\|{\Psi}\right\|^\mu=\sup_{k,q\in {\bf N}_0}\sup_{x\in I}\frac{\left|x^{k} (x^{-1}d/dx)^qx^{-\mu-1/2}\Psi(x)\right|}{(A_7+\sigma_1)^k a^{\ast 3}_k b^\ast_k (B_6+\rho_1)^q a^{\ast 2}_q b^{\ast 2}_q}\leq R \left\|h_{\mu} (\phi\#\psi)\right\|^\mu.
\end{eqnarray*}
Hence
$\Psi(x)\in H_{\mu,a^{\ast 3}_k b^{\ast}_k,A_7}^{a^{\ast 2}_q  b^{\ast 2}_q,B_6}$.\\

{\bf REMARK:} In Theorem 3.3 we may also choose $\phi,\psi \in H_{\mu,a_k,A}(or \widetilde{H}_{\mu}^{b_q,B})$  then p.d.o $h_{\mu,a}$ on Hankel convolution 
$(h_{\mu,a}(\phi\#\psi)) \in \widetilde{H}_{\mu,a^{\ast 2}_k,A'_7}(or H_{\mu}^{b^{\ast 2}_q,B'_6})$, where $ A'_7=(a^\ast_0/a^\ast_1)D+B_1$ and $B'_6=( b^\ast_0 / b^\ast_1 ) C+H_2^6 B^2$.
Similarly we may define $ \phi\#\psi, h_\mu(\phi\#\psi)$ for $\phi \in H_{\mu,a_k,A}$ and $\psi \in \widetilde{H}_\mu^{b_q,B}$ or $\psi \in H_{\mu,a_k,A}^{b_q,B}$ and study the p.d.o. on $\phi\#\psi$.\\

{\bf Acknowledgment}\\

This work has been supported by University Grants Commission, Govt. of India, under Grant No. F. No.  34-145$\backslash$2008  (SR).

\end{document}